%% file: conference_101719.tex
\documentclass[conference]{IEEEtran}
\IEEEoverridecommandlockouts
\usepackage{cite}
\usepackage{amsmath,amssymb,amsfonts}
\usepackage{algorithmic}
\usepackage{graphicx}
\usepackage{textcomp}
\usepackage{xcolor}
\usepackage{tikz}
\usepackage[version=4]{mhchem}
\usepackage{booktabs}
\usepackage{multirow}
\usepackage{siunitx}
\usepackage{eurosym}
\usepackage{mhchem}
\usepackage{amssymb}
\usepackage{tikz}
\usepackage{pgfplots}
\usepackage{multirow}
\usepackage{graphicx}
\usepackage{hyperref}
\usepackage{eurosym}
\def\BibTeX{{\rm B\kern-.05em{\sc i\kern-.025em b}\kern-.08em
    T\kern-.1667em\lower.7ex\hbox{E}\kern-.125emX}}
\usepackage{url}

\newcommand{\figspace}{\vspace{-0.5cm}}

\begin{document}

\title{Electrolyzer Scheduling for Nordic FCR Services\\
}

\author{%
\IEEEauthorblockN{Marco Saretta, Enrica Raheli, and Jalal Kazempour}
\IEEEauthorblockA{Department of Wind and Energy Systems, Technical University of Denmark, Kgs. Lyngby, Denmark\\ 
marco.saretta@outlook.com, $\{$enrah, jalal$\}$@dtu.dk 
}
}
\IEEEaftertitletext{\vspace{-0.8\baselineskip}}
\maketitle
\thispagestyle{plain}
\pagestyle{plain}
\begin{abstract}
The cost competitiveness of green hydrogen production via electrolysis presents a significant challenge for its large-scale adoption. One potential solution to make electrolyzers profitable is to diversify their products and participate in various markets, generating additional revenue streams.
Electrolyzers can be utilized as flexible loads and participate in various frequency-supporting ancillary service markets by adjusting their operating set points.
This paper develops a mixed-integer linear model, deriving an optimal scheduling strategy for an electrolyzer providing Frequency Containment Reserve (FCR) services in the Nordic synchronous region. Depending on the hydrogen price and demand, results show that the provision of various FCR services, particularly those for critical frequency conditions (FCR-D), could significantly increase the profit of the electrolyzer.  
\end{abstract}

\linespread{0.89}
\input{Sections/01_Introduction}
\input{Sections/02_Preliminaries}

\input{Sections/03_Model}
\input{Sections/04_Results}
\input{Sections/05_Conclusion}
\section*{Acknowledgment}
The authors would like to thank Thomas Dalgas Fechtenburg (Energinet) for our discussions on the setup and requirements for an electrolyzer to provide FCR services, Roar Hestbek Nicolaisen and Andreas Thingvad (Hybrid Greentech) for the support in the conceptualization and formulation of the model, and finally, Edoardo Simioni (\O rsted) for examining the MSc thesis project and providing constructive feedback.


\bibliography{bibliography.bib} 
\bibliographystyle{ieeetr}
\end{document}

%% file: Sections/01_Introduction.tex
\section{Introduction}

The production of renewable hydrogen through electrolysis is widely acknowledged as a crucial step in the green transition, enabling decarbonization of hard-to-abate sectors, such as industry and heavy transport. To support the large-scale development of electrolyzers, several countries in Europe and globally have released national hydrogen strategies. For example, in 2021 the Danish government released a strategy for the national development of Power-to-X, with a goal to construct 4 to 6 GW of electrolysis capacity by 2030 \cite{PtX_strategy_DK}. However, there are numerous challenges to scaling up this technology, including the cost competitiveness of the electrolysis-based hydrogen production \cite{GUERRA20192425}. This requires the establishment of new business models by diversification of the products \cite{newbusinessmodels}.

Electrolyzers are flexible assets that can rapidly change their power consumption level within their operating range with ramp rates around 20\% of the nominal power per second \cite{VARELA20219303,BUTTLER20182440}. This makes them eligible to produce various frequency-supporting ancillary services, providing an additional promising revenue stream 
\cite{Mancarella}. Examples of potential ancillary services that electrolyzers can produce are Frequency Containment Reserve (FCR) as a primary reserve, automatic Frequency Restoration Reserve (aFRR) as a secondary reserve, and manual Frequency Restoration Reserve (mFRR) as a tertiary reserve. The technical feasibility of electrolyzers for providing various services is investigated in \cite{Tuinema} and \cite{ALLIDIERES20199690}. The economic feasibility of providing grid services is analyzed in \cite{France} and \cite{Germany_offshore} for the French and German context, respectively. In \cite{ZHENG2023}, a scheduling model for an electrolyzer in Western Denmark (DK1) participating in the day-ahead, balancing, and reserve markets is proposed, showing that offering FCR and aFRR services significantly increases the profit. In a similar direction but for batteries, \cite{thingvad2022economic} develops a business model by selling FCR services in Eastern Denmark (DK2). All these studies show that the extent of increased profit by selling ancillary services depends significantly on the location of the electrolyzer due to different market products, prices, and eligibility requirements. 

This paper develops a scheduling model for an electrolyzer located in DK2, which is part of the Nordic synchronous region. Compared to the Continental Europe region including DK1, the power system in the Nordic region is smaller in scale and capacity, with a higher penetration rate of renewables, and thereby lower inertia. For that, there are three sub-categories of FCR services in the Nordic region designed for different ranges of frequency deviation, including FCR-N (for normal operations) and FCR-D Up/Down (for operations under disturbance with critically low/high frequency). 
The main contributions of this paper are twofold. First, we develop a mixed-integer linear model for scheduling electrolyzers, aiming to maximize their profit by selling hydrogen as well as FCR-N and FCR-D Up/Down services. 
Second, we provide a quantitative assessment to evaluate to what extent an electrolyzer located in DK2 earns more by providing FCR services, in comparison to  a case that solely produces hydrogen. 

The remaining of the paper is organized as follows. Section \ref{preliminaries} provides an introduction to the Nordic FCR markets. Section \ref{model} presents the proposed optimization model. Section \ref{results} provides numerical scheduling results and an economic assessment. Finally, Section \ref{conclusion} concludes the paper.

\textit{Notation:} By $\lambda^{\text{(.)}}_t$, we refer to forecast for the volume-weighted average price  $\lambda^{\text{FCR-N}}_t$ for FCR-N,  $\lambda^{\text{FCR-D}\uparrow}_t$ for FCR-D Up, and  $\lambda^{\text{FCR-D}\downarrow}_t$ for FCR-D Down, all in hour $t$. Similarly, let $r^{\text{(.)}}_t$ denote the quantity bids $r^{\text{FCR-N}}_t$, $r^{\text{FCR-D}\uparrow}_t$, and $r^{\text{FCR-D}\downarrow}_t$ to be submitted to the corresponding markets. 




%% file: Sections/02_Preliminaries.tex
\section{Preliminaries: Nordic FCR markets}
\label{preliminaries}
\subsection{General overview}
The Transmission System Operator (TSO) is the organization
in charge on a national scale for the secure operation of the power grid. TSOs within synchronous areas share responsibility for real-time balance between supply and demand to maintain the grid frequency close to the nominal value, e.g., 50 Hz in Europe. Ancillary services are the measures adopted by TSOs to ensure grid stability. For that, TSOs procure reserves for ancillary services in advance, and activate them in the real-time operation if necessary. 

For completeness, Table~\ref{tab:fcr nomen} provides a  nomenclature for frequency-supporting ancillary services in the Nordic and Continental Europe synchronous regions, although the focus of this paper is the FCR services in the Nordic region.


\begin{table}[]
\figspace
\centering
\caption{Nomenclature for ancillary services  in the continental Europe and Nordic synchronous regions \cite{energinet2012ancillary}}
\resizebox{\columnwidth}{!}{%
\begin{tabular}{@{}ccc@{}}
\toprule
Function                                                                                             & Continental Europe (including DK1)                                                                             & Nordics (including DK2)                                                                                       \\ \midrule
\multirow{3}{*}{\begin{tabular}[c]{@{}c@{}}Frequency stabilization\\ (primary reserve)\end{tabular}} & \multirow{3}{*}{\begin{tabular}[c]{@{}c@{}}Frequency Containment Reserve\\ (FCR)\end{tabular}} & \begin{tabular}[c]{@{}c@{}}FCR for Normal Operations \\ (FCR-N)\end{tabular}                   \\ \cmidrule(l){3-3} 
                                                                                                     &                                                                                                & \begin{tabular}[c]{@{}c@{}}FCR for Disturbances - Up regulation\\ (FCR-D Up)\end{tabular}     \\ \cmidrule(l){3-3} 
                                                                                                     &                                                                                                & \begin{tabular}[c]{@{}c@{}}FCR for Disturbances - Down regulation\\ (FCR-D Down)\end{tabular} \\ \midrule
\begin{tabular}[c]{@{}c@{}}Frequency recovery\\ (secondary reserve)\end{tabular}                     & \begin{tabular}[c]{@{}c@{}}automatic Frequency Restoration Reserve\\ (aFRR)\end{tabular}       & \begin{tabular}[c]{@{}c@{}}automatic Frequency Restoration Reserve\\ (aFRR)\end{tabular}      \\ \midrule
\begin{tabular}[c]{@{}c@{}}Balance adjustment\\ (tertiary reserve)\end{tabular}                      & \begin{tabular}[c]{@{}c@{}}manual Frequency Restoration Reserve\\ (mFRR)\end{tabular}          & \begin{tabular}[c]{@{}c@{}}manual Frequency Restoration Reserve\\ (mFRR)\end{tabular}         \\ \bottomrule
\end{tabular}%
}
\label{tab:fcr nomen}
\figspace
\end{table}
\subsection{Market structure}
\begin{figure}[b]
\figspace
    \centering
\tikzset{every picture/.style={line width=0.75pt}} 
\resizebox{0.45\textwidth}{!}{%
\begin{tikzpicture}[x=0.75pt,y=0.75pt,yscale=-1,xscale=1]

\draw [color={rgb, 255:red, 0; green, 0; blue, 0 }  ,draw opacity=1 ][fill={rgb, 255:red, 245; green, 166; blue, 35 }  ,fill opacity=0.13 ][line width=1.5]    (71,161) -- (635.25,161) (81,157) -- (81,165)(91,157) -- (91,165)(101,157) -- (101,165)(111,157) -- (111,165)(121,157) -- (121,165)(131,157) -- (131,165)(141,157) -- (141,165)(151,157) -- (151,165)(161,157) -- (161,165)(171,157) -- (171,165)(181,157) -- (181,165)(191,157) -- (191,165)(201,157) -- (201,165)(211,157) -- (211,165)(221,157) -- (221,165)(231,157) -- (231,165)(241,157) -- (241,165)(251,157) -- (251,165)(261,157) -- (261,165)(271,157) -- (271,165)(281,157) -- (281,165)(291,157) -- (291,165)(301,157) -- (301,165)(311,157) -- (311,165)(321,157) -- (321,165)(331,157) -- (331,165)(341,157) -- (341,165)(351,157) -- (351,165)(361,157) -- (361,165)(371,157) -- (371,165)(381,157) -- (381,165)(391,157) -- (391,165)(401,157) -- (401,165)(411,157) -- (411,165)(421,157) -- (421,165)(431,157) -- (431,165)(441,157) -- (441,165)(451,157) -- (451,165)(461,157) -- (461,165)(471,157) -- (471,165)(481,157) -- (481,165)(491,157) -- (491,165)(501,157) -- (501,165)(511,157) -- (511,165)(521,157) -- (521,165)(531,157) -- (531,165)(541,157) -- (541,165)(551,157) -- (551,165)(561,157) -- (561,165)(571,157) -- (571,165)(581,157) -- (581,165)(591,157) -- (591,165)(601,157) -- (601,165)(611,157) -- (611,165)(621,157) -- (621,165)(631,157) -- (631,165) ;
\draw [shift={(638.25,161)}, rotate = 180] [color={rgb, 255:red, 0; green, 0; blue, 0 }  ,draw opacity=1 ][line width=1.5]    (14.21,-4.28) .. controls (9.04,-1.82) and (4.3,-0.39) .. (0,0) .. controls (4.3,0.39) and (9.04,1.82) .. (14.21,4.28)   ;
\draw  [color={rgb, 255:red, 74; green, 144; blue, 226 }  ,draw opacity=1 ][fill={rgb, 255:red, 74; green, 144; blue, 226 }  ,fill opacity=0.13 ][line width=1.5]  (71,161) -- (310.75,161) -- (310.75,201) -- (71,201) -- cycle ;
\draw [color={rgb, 255:red, 74; green, 144; blue, 226 }  ,draw opacity=1 ][fill={rgb, 255:red, 74; green, 144; blue, 226 }  ,fill opacity=0.13 ][line width=1.5]  [dash pattern={on 1.69pt off 2.76pt}]  (210.72,90.85) -- (210.72,136.85) ;
\draw [color={rgb, 255:red, 74; green, 144; blue, 226 }  ,draw opacity=1 ][fill={rgb, 255:red, 74; green, 144; blue, 226 }  ,fill opacity=0.13 ][line width=1.5]  [dash pattern={on 1.69pt off 2.76pt}]  (221.5,56.95) -- (221.5,113.85) ;
\draw [color={rgb, 255:red, 74; green, 144; blue, 226 }  ,draw opacity=1 ][fill={rgb, 255:red, 74; green, 144; blue, 226 }  ,fill opacity=0.13 ][line width=1.5]  [dash pattern={on 1.69pt off 2.76pt}]  (221.5,135.85) -- (221.5,155.25) ;
\draw  [color={rgb, 255:red, 74; green, 144; blue, 226 }  ,draw opacity=1 ][fill={rgb, 255:red, 74; green, 144; blue, 226 }  ,fill opacity=0.13 ][line width=1.5]  (213.72,118.74) .. controls (213.72,117.02) and (215.11,115.63) .. (216.83,115.63) -- (226.17,115.63) .. controls (227.89,115.63) and (229.28,117.02) .. (229.28,118.74) -- (229.28,128.11) .. controls (229.28,129.83) and (227.89,131.23) .. (226.17,131.23) -- (216.83,131.23) .. controls (215.11,131.23) and (213.72,129.83) .. (213.72,128.11) -- cycle ;

\draw  [color={rgb, 255:red, 74; green, 144; blue, 226 }  ,draw opacity=1 ][fill={rgb, 255:red, 74; green, 144; blue, 226 }  ,fill opacity=0.13 ][line width=1.5]  (157,23) .. controls (157,18.58) and (160.58,15) .. (165,15) -- (278,15) .. controls (282.42,15) and (286,18.58) .. (286,23) -- (286,47) .. controls (286,51.42) and (282.42,55) .. (278,55) -- (165,55) .. controls (160.58,55) and (157,51.42) .. (157,47) -- cycle ;

\draw [color={rgb, 255:red, 248; green, 180; blue, 28 }  ,draw opacity=1 ][fill={rgb, 255:red, 245; green, 166; blue, 35 }  ,fill opacity=0.13 ][line width=1.5]  [dash pattern={on 1.69pt off 2.76pt}]  (490.57,54.45) -- (490.57,135.25) ;
\draw  [color={rgb, 255:red, 248; green, 180; blue, 28 }  ,draw opacity=1 ][fill={rgb, 255:red, 245; green, 166; blue, 35 }  ,fill opacity=0.13 ][line width=1.5]  (482.84,141.36) .. controls (482.84,139.64) and (484.24,138.25) .. (485.95,138.25) -- (495.3,138.25) .. controls (497.01,138.25) and (498.41,139.64) .. (498.41,141.36) -- (498.41,150.74) .. controls (498.41,152.46) and (497.01,153.85) .. (495.3,153.85) -- (485.95,153.85) .. controls (484.24,153.85) and (482.84,152.46) .. (482.84,150.74) -- cycle ;

\draw  [color={rgb, 255:red, 248; green, 180; blue, 28 }  ,draw opacity=1 ][fill={rgb, 255:red, 245; green, 166; blue, 35 }  ,fill opacity=0.13 ][line width=1.5]  (412,19.29) .. controls (412,14.71) and (415.71,11) .. (420.29,11) -- (560.96,11) .. controls (565.54,11) and (569.25,14.71) .. (569.25,19.29) -- (569.25,44.16) .. controls (569.25,48.74) and (565.54,52.45) .. (560.96,52.45) -- (420.29,52.45) .. controls (415.71,52.45) and (412,48.74) .. (412,44.16) -- cycle ;

\draw  [color={rgb, 255:red, 248; green, 180; blue, 28 }  ,draw opacity=1 ][fill={rgb, 255:red, 245; green, 166; blue, 35 }  ,fill opacity=0.13 ][line width=1.5]  (310.75,161) -- (551,161) -- (551,201) -- (310.75,201) -- cycle ;
\draw [color={rgb, 255:red, 74; green, 144; blue, 226 }  ,draw opacity=1 ][fill={rgb, 255:red, 74; green, 144; blue, 226 }  ,fill opacity=0.13 ][line width=1.5]  [dash pattern={on 1.69pt off 2.76pt}]  (210.72,90.85) -- (127,90.95) ;
\draw  [color={rgb, 255:red, 74; green, 144; blue, 226 }  ,draw opacity=1 ][fill={rgb, 255:red, 74; green, 144; blue, 226 }  ,fill opacity=0.13 ][line width=1.5]  (202.94,142.36) .. controls (202.94,140.64) and (204.34,139.25) .. (206.05,139.25) -- (215.4,139.25) .. controls (217.11,139.25) and (218.51,140.64) .. (218.51,142.36) -- (218.51,151.74) .. controls (218.51,153.46) and (217.11,154.85) .. (215.4,154.85) -- (206.05,154.85) .. controls (204.34,154.85) and (202.94,153.46) .. (202.94,151.74) -- cycle ;

\draw  [color={rgb, 255:red, 74; green, 144; blue, 226 }  ,draw opacity=1 ][fill={rgb, 255:red, 74; green, 144; blue, 226 }  ,fill opacity=0.13 ][line width=1.5]  (7,78.5) .. controls (7,74.08) and (10.58,70.5) .. (15,70.5) -- (116,70.5) .. controls (120.42,70.5) and (124,74.08) .. (124,78.5) -- (124,102.5) .. controls (124,106.92) and (120.42,110.5) .. (116,110.5) -- (15,110.5) .. controls (10.58,110.5) and (7,106.92) .. (7,102.5) -- cycle ;

\draw  [color={rgb, 255:red, 248; green, 180; blue, 28 }  ,draw opacity=1 ][fill={rgb, 255:red, 245; green, 166; blue, 35 }  ,fill opacity=0.13 ][line width=1.5]  (299.5,76.5) .. controls (299.5,72.08) and (303.08,68.5) .. (307.5,68.5) -- (408.5,68.5) .. controls (412.92,68.5) and (416.5,72.08) .. (416.5,76.5) -- (416.5,100.5) .. controls (416.5,104.92) and (412.92,108.5) .. (408.5,108.5) -- (307.5,108.5) .. controls (303.08,108.5) and (299.5,104.92) .. (299.5,100.5) -- cycle ;

\draw [color={rgb, 255:red, 248; green, 180; blue, 28 }  ,draw opacity=1 ][fill={rgb, 255:red, 245; green, 166; blue, 35 }  ,fill opacity=0.13 ][line width=1.5]  [dash pattern={on 1.69pt off 2.76pt}]  (450.72,89.85) -- (450.72,135.85) ;
\draw [color={rgb, 255:red, 248; green, 180; blue, 28 }  ,draw opacity=1 ][fill={rgb, 255:red, 245; green, 166; blue, 35 }  ,fill opacity=0.13 ][line width=1.5]  [dash pattern={on 1.69pt off 2.76pt}]  (450.72,89.85) -- (431,89.87) -- (417,89.87) ;
\draw  [color={rgb, 255:red, 248; green, 180; blue, 28 }  ,draw opacity=1 ][fill={rgb, 255:red, 245; green, 166; blue, 35 }  ,fill opacity=0.13 ][line width=1.5]  (442.94,141.36) .. controls (442.94,139.64) and (444.34,138.25) .. (446.05,138.25) -- (455.4,138.25) .. controls (457.11,138.25) and (458.51,139.64) .. (458.51,141.36) -- (458.51,150.74) .. controls (458.51,152.46) and (457.11,153.85) .. (455.4,153.85) -- (446.05,153.85) .. controls (444.34,153.85) and (442.94,152.46) .. (442.94,150.74) -- cycle ;

\draw  [color={rgb, 255:red, 57; green, 211; blue, 33 }  ,draw opacity=1 ][fill={rgb, 255:red, 149; green, 233; blue, 134 }  ,fill opacity=0.23 ][line width=1.5]  (551,161) -- (611.25,161) -- (611.25,201) -- (551,201) -- cycle ;

\draw (175.44,172.5) node [anchor=north west][inner sep=0.75pt]   [align=left] {$\rm{D}$-$2$};
\draw (410.32,172.5) node [anchor=north west][inner sep=0.75pt]   [align=left] {$\rm{D}$-$1$};
\draw (574,173) node [anchor=north west][inner sep=0.75pt]   [align=left] {$\rm{D}$};
\draw (214,117.43) node [anchor=north west][inner sep=0.75pt]  [font=\footnotesize] [align=left] {15};
\draw (203.22,141) node [anchor=north west][inner sep=0.75pt]  [font=\footnotesize] [align=left] {14};
\draw (11.5,74) node [anchor=north west][inner sep=0.75pt]  [font=\small] [align=left] {\begin{minipage}[lt]{75.14pt}\setlength\topsep{0pt}
\begin{center}
Published spot \\prices for day $\rm{D}$-$1$
\end{center}

\end{minipage}};
\draw (310,72) node [anchor=north west][inner sep=0.75pt]  [font=\small] [align=left] {\begin{minipage}[lt]{66.98pt}\setlength\topsep{0pt}
\begin{center}
Published spot \\prices for day $\rm{D}$
\end{center}

\end{minipage}};
\draw (170,18.5) node [anchor=north west][inner sep=0.75pt]  [font=\small] [align=left] {FCR first auction\\closure for day $\rm{D}$};
\draw (443.22,140) node [anchor=north west][inner sep=0.75pt]  [font=\footnotesize] [align=left] {14};
\draw (483.13,140) node [anchor=north west][inner sep=0.75pt]  [font=\footnotesize] [align=left] {18};
\draw (426.63,15.22) node [anchor=north west][inner sep=0.75pt]  [font=\small] [align=left] {\begin{minipage}[lt]{88.4pt}\setlength\topsep{0pt}
\begin{center}
FCR second auction \\closure for day $\rm{D}$
\end{center}

\end{minipage}};
\draw (605,136) node [anchor=north west][inner sep=0.75pt]   [align=left] {[hour]};

\end{tikzpicture}
}
\vspace{-0.5mm}
    \caption{Timeline for FCR and spot markets in the Nordic region. There are two auctions for the FCR services, one before and one after the spot market.}
    \label{fig:timeline_auction}
    \figspace
\end{figure}
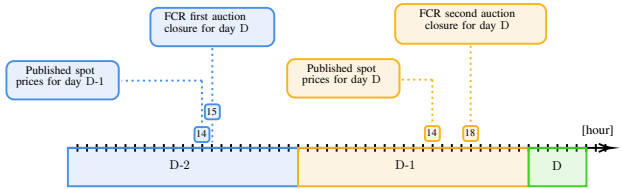

The Nordic obligations indicate the reserve requirements that must be collectively secured in every hour among the Nordic TSOs in a proportional share for different services, as reported in Table~\ref{tab:nordic_obligations_table}.
Note that StatNett, FinGrid, Svenska Kraftnat, and Energinet are national TSOs in Norway, Finland, Sweden, and Denmark, respectively. 
The Danish TSO, Energinet, has a comparatively lower share due to congestion and technical limitations of the DK2-Sweden connection cable. The Nordic obligations for any hour of the day $\rm{D}$ are contracted via two separate auctions, both pay-as-bid structured, on $\rm{D}$-$2$ and $\rm{D}$-$1$ prior to the delivery day, as shown in Figure~\ref{fig:timeline_auction}. Approximately, 80\% of each FCR service is contracted on the $\rm{D}$-$2$ auction, and the remaining in $\rm{D}$-$1$.  

During the daily FCR auctions, Energinet and Svenska Kraftnat jointly procure their share of reserves, hence Danish FCR providers can potentially meet the full Swedish demand for FCR-N and FCR-D services. However, the maximum amount of FCR from a single unit is limited to 100 MW \cite{energinet2012ancillary} to avoid a significant loss of FCR in case of a unit failure.
\subsection{FCR delivery and payment structure}
The provision of FCR services entails two distinct stages, namely reserve contraction and activation.

Reserve contraction occurs during the $\rm{D}$-$2$ or $\rm{D}$-$1$ auction, wherein the availability of the reserve noted in the FCR bid is approved by the TSO. Recall that both auctions are based on a pay-as-bid scheme. Compensation for the FCR service in hour $t$ is based on the reserve quantity $r^{\text{(.)}}_t$ (MW) and the submitted bid price $\lambda^{\text{(.)}}_t$ (\euro/MW), resulting a revenue, the so-called reserve payment. The Nordic TSOs do not currently disclose information about the last accepted bid in the auctions. The only public information is the hourly volume-weighted average bid price for each service once the auction is closed.

 \begin{table}[t]
\figspace
\centering
\caption{Nordic obligations for FCR services in  2023 \cite{energinet2012ancillary}}
\resizebox{0.45\textwidth}{!}{%
\begin{tabular}{lcccc}
\toprule
 & Share   & FCR-N    & FCR-D Up &  FCR-D Down              \\ 
& [\%] & [MW] & [MW] &  [MW]                  \\ \midrule
StatNett                                 & 39                           & 234                          & 564                           & 546                                         \\
FinGrid                                  &  20                           & 120                          & 290                           & 280                                         \\
Svenska Kraftnat                         & 38.3                         & 230                          & 555                           & 536                                       \\
Energinet                                &  2.7                          & 17                           & 41                            & 38                                        \\ \midrule
Nordic obligations                       &  100                          & 600                          & 1450                          & 1400                                        \\ \bottomrule
\end{tabular}%
}
\label{tab:nordic_obligations_table}
\figspace
\end{table}



The activation payment is linked to the real-time operation, where the FCR provider must activate the reserve according to the frequency level $f$ in Hz at any instant within the hour declared in the bid. The real-time reserve activation at any instant in hour $t$ is equal to the product of the amount of the contracted reserve $r^{\text{(.)}}_t$ and the normalized instantaneous response $y^{\text{(.)}}$, defined below for FCR-N, FCR-D Up, and FCR-D Down, respectively:
\vspace{-0.3mm}
\begin{subequations} \label{eq:FCR}
\begin{equation}
\resizebox{0.275\textwidth}{!}{
$\displaystyle
y^{\text{FCR-N}}=\left\{\begin{aligned}
-1, & \quad \text {if} \quad f < \SI{49.9}{} \\
\frac{f - \SI{50}{}}{\SI{0.1}{}}, & \quad \text {if} \quad \SI{49.9}{} \leq f \leq \SI{50.1}{} \\
+1, & \quad \text {if} \quad f > \SI{50.1}{}
\end{aligned}\right.
$
\label{eq:FCRN_eq}
}
\end{equation}
\begin{equation}
\resizebox{0.275\textwidth}{!}{
$\displaystyle
y^{\text{FCR-D}\uparrow}=\left\{\begin{aligned}
-1, & \quad \text {if} \quad f < \SI{49.5}{} \\
\frac{f - \SI{49.9}{}}{\SI{0.4}{}}, & \quad \text {if} \quad \SI{49.5}{} \leq f \leq \SI{49.9}{} \\
0, & \quad \text {if} \quad f > \SI{49.9}{}
\end{aligned}\right.
$
\label{eq:FCRD_up_equation}
}
\end{equation}
\begin{equation}
\resizebox{0.275\textwidth}{!}{
$\displaystyle
y^{\text{FCR-D}\downarrow}=\left\{\begin{aligned}
0, & \quad \text {if} \quad f < \SI{50.1}{} \\
\frac{f - \SI{50.1}{}}{\SI{0.4}{}}, & \quad \text {if} \quad \SI{50.1}{} \leq f \leq \SI{50.5}{} \\
+1, & \quad \text {if} \quad f > \SI{50.5}{}.
\end{aligned}\right.
$
\label{eq:fcrd_down_equation}
}
\end{equation}
\end{subequations} 
\vspace{-0.3mm}
The payment for activated quantity $r^{\text{(.)}}_t y^{\text{(.)}}$ at any instant within hour $t$ is based on the balancing price in the corresponding hour. The settlement typically occurs by the TSO within one week after the service delivery.

\vspace{-0.6mm}
\subsection{Electrolyzer eligibility assessment} 
To qualify for FCR service provision, the Nordic TSOs have established pre-qualification requirements in terms of response time. For FCR-D Up/Down, the electrolyzer must be capable to respond for the half reserve within 5 seconds, and the full reserve within 30 seconds. For FCR-N, the full reserve must be activated within 150 seconds. To determine if an alkaline electrolyzer is eligible, its ramp-rate compliance needs to be assessed. Manufacturers usually do not disclose ramp rates information. However, an estimation of a ramp-rate around 20\% of the nominal capacity per second makes an alkaline electrolyzer eligible for FCR-N and FCR-D Up/Down provision \cite{VARELA20219303}.

%% file: Sections/03_Model.tex
\section{Proposed Optimization Model}
\label{model}
\begin{figure}[t]
\centering
\begin{tikzpicture}
\node[inner sep=0pt] (grid) at (0,0)
    {\includegraphics[width=.03\textwidth]{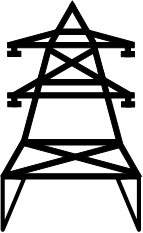}};
\node[inner sep=0pt] (electrolyzer) at (1.9,0)
   [draw,thick,align=center,minimum width=1.8cm,minimum height=0.9cm, text width = 1.5 cm, label=below:{ $z^{\mathrm{on}}_{t},z^{\mathrm{sb}}_{t},z^{\mathrm{su}}_{t}$}] {\footnotesize Electrolyzer system};
\node[inner sep=0pt] (FCR) at (1.9,1.5)
   [draw,thick,align=center,minimum width=1.8cm,minimum height=0.9cm, text width = 1.5 cm] {\footnotesize FCR auctions};
\node[inner sep=0pt] (compressor) at (3.8,0)
    {\includegraphics[width=.04\textwidth]{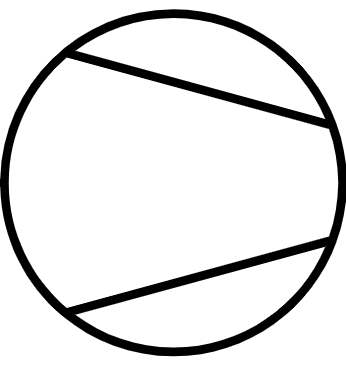}};
\node[inner sep=0pt] (compressor_above) at (3.8,0.5) [label=above:{\footnotesize Compressor}]{};
\node[inner sep=0pt,fill,circle,minimum size=0.1 cm] (after_compr) at (4.5,0) {};
\node[inner sep=0pt] (storage) at (5.35,0.5) [label=above:{\footnotesize  Storage}]
    {\includegraphics[width=.04\textwidth]{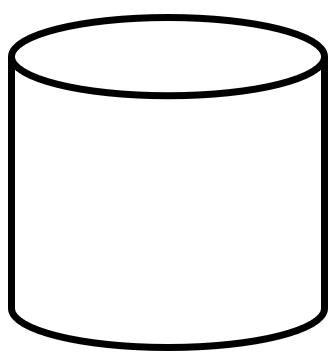}};
\node[inner sep=0pt,fill,circle,minimum size=0.1 cm] (before_demand) at (6.2,0) {};
\node[inner sep=0pt] (demand) at (7.8,0) [align=center,minimum width=1.5cm,minimum height=0.9cm, text width = 1.5 cm] {\footnotesize Hydrogen demand};
\draw[->,thick] (grid.east) -- (electrolyzer.west)
    node[midway,above] {$p_t$};
\draw[->,thick] (electrolyzer.north) -- (FCR.south)
    node[midway,right] {$r_t^{\text{(.)}}$};
\draw[->,thick] (electrolyzer.east) -- (compressor.west)
    node[midway,above] {$h_t^{\rm{p}}$};
\draw[->,thick] (compressor.east) -- (after_compr.west)
    node[midway,above] {};
\draw[->,thick] (before_demand.east) -- (demand.west)
    node[midway,above] {$d_t$};

\node[inner sep=0pt,fill,circle,minimum size=0 cm] (A) at (4.5,0.5) {};
\node[inner sep=0pt,fill,circle,minimum size=0 cm] (B) at (4.5,-0.5) {};
\node[inner sep=0pt,fill,circle,minimum size=0 cm] (C) at (6.2,0.5) {};
\node[inner sep=0pt,fill,circle,minimum size=0 cm] (D) at (6.2,-0.5) {};
\draw[-,thick] (after_compr.north) -- (A.south) node[midway,above] {};
\draw[->,thick] (A.east) -- (storage.west) node[midway,above] {};
\draw[-,thick] (storage.east) -- (C.west) node[midway,above] {};
\draw[->,thick] (C.south) -- (before_demand.north) node[midway,above] {};
\draw[-,thick] (after_compr.south) -- (B.north) node[midway,above] {};
\draw[-,thick] (B.east) -- (D.west) node[midway,above] {};
\draw[->,thick] (D.north) -- (before_demand.south) node[midway,above] {};
\end{tikzpicture}
\figspace
\caption{Power system design of an electrolyzer system and its auxiliary assets, providing hydrogen and FCR services. The power source is the grid.}
\label{fig:power_diagram}
\figspace
\end{figure}
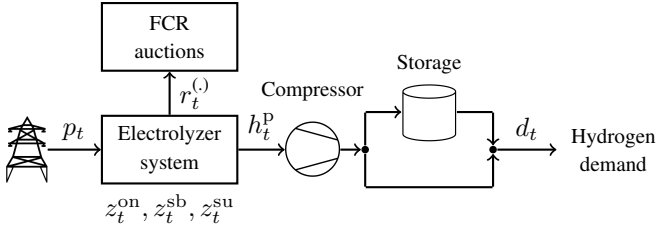

This paper focuses on a typical alkaline electrolyzer providing hydrogen and FCR services, while purchasing power from the the grid --- we do not consider any local renewable power supply. This allows for a constant baseline during the FCR scheduling and an adjustment in power consumption according to \eqref{eq:FCR} when an activation is required. Figure~\ref{fig:power_diagram} illustrates the electrolyzer system and its auxiliary assets including hydrogen compressor and storage.  
The  hydrogen produced is compressed and delivered to an off-taker (demand), meeting a weekly demand, sold at a fixed price (\euro/kg).

\subsection{When do we solve the optimization problem?}
Recall from Figure~\ref{fig:timeline_auction} that the electrolyzer has the opportunity to participate in two pay-as-bid auctions for FCR services, one being settled in day $\rm{D}$-$2$ and the other one in $\rm{D}$-$1$. The electrolyzer owner can solve our proposed optimization model for scheduling decision making at two distinct points of time:
\begin{enumerate}
    \item At any time before the first auction closure at hour 15:00 of day $\rm{D}$-$2$. In this case, the electrolyzer owner should forecast prices $\lambda^{\text{(.)}}_t$ for the first auction, as well as the hourly spot prices $\lambda_{t}^{\mathrm{spot}}$ whose true values will be realized in day $\rm{D}$-$1$. All these forecasted prices are treated as  input parameters to our optimization model.  By solving it, we will determine reserve quantity bids $r^{\text{(.)}}_t$ to be submitted to the first FCR auction. The price bids are the same as the forecasted prices. This optimization problem also gives quantity bids to be submitted to the spot market, i.e., $p_{t}$,  but they can be modified until noon of $\rm{D}$-$1$ by a re-optimization with fixed $r^{\text{(.)}}_t$ if updated spot price forecasts are available. 
    
    \item At any time before noon of day $\rm{D}$-$1$, i.e., the closure of the spot market, in case the electrolyzer could not sell FCR services in the first auction. This time, $\lambda^{\text{(.)}}_t$ are forecasted prices for the second auction, which are not necessarily identical to realized prices of the first auction. The optimization outcomes are quantity bids to be submitted to the spot market, i.e., $p_{t}$, and to the second FCR auction, i.e., $r^{\text{(.)}}_t$.  Note that we can also solve this optimization problem between hours 14 and 18 of day $\rm{D}$-$1$, but then the hourly power purchases $p_{t}$ are fixed based on the spot market outcomes, and therefore hourly power consumptions can no longer be changed unless by trading in the intra-day and subsequent markets, which is outside the scope of this paper. 
\end{enumerate}
\subsection{Mathematical formulation}
The proposed  model is formulated as \eqref{eq:2_objective}-\eqref{eq:var}. Lower-case symbols are used for variables, whereas upper-case or Greek symbols indicate parameters. The objective function maximizes the total profit over the set of hours $t \in \mathcal{T}$ as 
\begin{align}
\begin{split}
    & \underset{\bf{x}}{\rm{max}} \sum_{t \in \mathcal{T}} \bigg( d_{t}  \lambda^{{\ce{H2}}} +  r^{\text{FCR-N}}_t  \lambda^{\text{FCR-N}}_t + r^{\text{FCR-D}\uparrow}_t  \lambda^{\text{FCR-D}\uparrow}_t+\\ &  r^{\text{FCR-D}\downarrow}_t  \lambda^{\text{FCR-D}\downarrow}_t -p_{t}  \left( \lambda_{t}^{\mathrm{spot}} + \lambda^{\text{TSO}} + \lambda^{\text{DSO}} \right) -z^{\text{su}}  K^{\text{su}} \bigg),
   \label{eq:2_objective}
\end{split}
\end{align}
where vector $\bf{x}$ includes the set of variables, which will be defined later. The revenue streams are based on hydrogen sale $d_{t}$ at a constant price $\lambda^{{\ce{H2}}}$ and service sales $r^{\text{(.)}}_t$ at price $\lambda^{\text{(.)}}_t$. The cost incurs by purchasing hourly power $p_{t}$ at spot market price $\lambda_{t}^{\mathrm{spot}}$, marked up by the TSO tariff $\lambda^{\text{TSO}}$ as well as the tariff of the distribution system operator $\lambda^{\text{DSO}}$, if the electrolyzer is comparatively small and connected to a distribution grid. In addition, \eqref{eq:2_objective} accounts for the cold start-up cost of the electrolyzer, where the binary variable $z^{\text{su}}$ indicates the start-up at hour $t$, associated with the cost $K^{\text{su}}$ per start-up. Note that the activation payment is excluded\footnote{This is a mild assumption because (\textit{i}) FCR-N is a service being activated in both sides. Historically, the FCR-N activation is almost symmetrical over every week in 2022, and (\textit{ii}) the activation rate of FCR-D Up/Down services in the Nordic area was less than 1\% in 2022 \cite{thingvad2022economic}. The interested reader in FCR activation data in DK2 is referred to \cite{fingrid}.}.




The following set of constraints \eqref{eq:elec} models the physics and limitations of the electrolyzer and auxiliary assets including compressor and hydrogen storage. The power purchased from the spot market, i.e., $p_{t}$, supplies the electrolyzer's consumption  $p^{\rm{e}}_{t}$ and the compressor's consumption  $p^{\rm{c}}_{t}$:
\begin{subequations} 
\begin{align}
    p_{t} &= p^{\rm{e}}_{t} + p^{\rm{c}}_{t} & \forall~ & t \in \mathcal{T}. \label{eq:P2X_balance}
\end{align}

The electrolyzer is either on, or standby, or off, i.e.,
\label{eq:elec}
\begin{align}
z^{\text{on}}_{t} +z^{\text{sb}}_{t}  &\leq 1 
& \forall \ t \in \mathcal{T},
\label{eq:exclusivity_states}
\end{align}
including binary variables $z^{\text{on}}_{t}$ (if 1, the electrolyzer is on) and $z^{\text{sb}}_{t} $ (if 1, the electrolyzer is on the standby state). If on, the electrolyzer consumes power and
produces hydrogen. If standby, the electrolyzer does not produce hydrogen but consumes 1-5\% of the nominal power needed to keep the system warm and pressurized for quick activation \cite{VARELA20219303}. If both binary variables are zero, then the electrolyzer is off, neither consuming power nor producing hydrogen\footnote{In this formulation, we model three states (on, standby, off) with two binary variables only, instead of three, as it is prevalent in the literature. We hypothesize, depending on the solver used, this may reduce computational time, but a further investigation is required.}.

The power consumption $p^{\text{e}}_{t}$ of the electrolyzer defines the operational baseline, constrained by 
\begin{align}
     P^{\text{min}}  z^{\text{on}}_{t} + P^{\text{sb}}  z^{\text{sb}}_{t} \leq p^{\text{e}}_{t} &\leq P^{\text{max}}  z^{\text{on}}_{t} + P^{\text{sb}}  z^{\text{sb}}_{t}& \forall~ & t \in \mathcal{T}, \label{eq:elyzer_disp_up}
\end{align}
where the lower bound is $P^{\text{min}}$ and the upper bound is the capacity $P^{\text{max}}$ when the electrolyzer is on ($z^{\text{on}}_{t}$=1). If standby ($z^{\text{sb}}_{t}$=1), $p^{\text{e}}_{t}$ is set to be equal to the standby power $P^{\text{sb}}$.

Transition from off state in hour $t-1$ to on state in $t$ incurs the start-up cost due to the need to reach the desired pressure and temperature levels. For that, \eqref{3states_cold} sets the binary variable $z^{\text{su}}_{t}$ to be 1 during such a transition, otherwise it is 0:
\begin{align}
     z^{\text{su}}_{t} &\geq (z^{\text{on}}_{t} - z^{\text{on}}_{t-1}) + (z^{\text{sb}}_{t} -z^{\text{sb}}_{t-1}) & \forall \ t \in \mathcal{T}.
     \label{3states_cold}
\end{align}

The power-to-hydrogen conversion efficiency of an alkaline electrolyzer is not constant over the operating range. To model the non-linear dependency between power consumption and hydrogen production, a piece-wise linearization is introduced as proposed in \cite{baumhof2023optimization}. For each linearization segment $s \in \mathcal{S}$, the hydrogen production $h^{\text{p}}_{t}$ is formulated as a linear function of the power consumption $\widehat{p}^{\text{e}}_{t, s}$:
\begin{align}
    h^{\text{p}}_{t} & = \sum_{s \in \mathcal{S}} \left( A_s  \widehat{p}^{\text{e}}_{t, s} + B_s \widehat{z}_{t, s} \right) & \forall \ t \in \mathcal{T},
\end{align}
where the coefficients $A_s$ and $B_s$ represent the slope and intercept for each linear segment, whereas the binary variable 
$\widehat{z}_{t, s}$, if one, indicates segment $s$ is active in hour $t$. 

The electrolyzer produces hydrogen in the on state with one segment active only in hour $t$, as enforced by
\begin{align}
    \sum_{s \in \mathcal{S}}&  \widehat{z}_{t, s} = z^{\text{on}}_t  & \forall \ t \in \mathcal{T} \label{piecewise_bin}.
\end{align}

The power consumption $\widehat{p}^{\text{e}}_{t, s}$ for each segment $s$ is constrained by
\begin{align}
   \underline{P}_s \widehat{z}_{t, s}&  \leq \widehat{p}^{\text{e}}_{t, s} \leq \overline{P}_s \widehat{z}_{t, s} & \forall \ t \in \mathcal{T}, \ \forall \ s \in \mathcal{S}, 
\end{align}
where $\underline{P}_s $ and $\overline{P}_s $  represent lower and upper bounds. The power consumption $p_t^{\text{e}}$ is then calculated as
\begin{align}
   p_t^{\text{e}} & = P^{\text{sb}}  z^{\text{sb}}_t + \sum_{s \in \mathcal{S}}  \widehat{p}^{\text{e}}_{t, s}   & \forall \ t \in \mathcal{T}.
\end{align}

The hydrogen production of the electrolyzer goes to the compressor to be further pressurized, and then is either stored or is directly injected to tube trailers, representing the demand. The compressor power consumption $p^{\rm{c}}_{t}$ is a function of the hydrogen production $h_t^{\mathrm{p}}$ of the electrolyzer as
\begin{align}
    p^{\rm{c}}_{t} &= K^{\rm{c}} h_t^{\mathrm{p}}   & \forall~ & t \in \mathcal{T}, \label{eq:compr}
\end{align}
where  $K^{\rm{c}}$ gives the energy required to compress 1 kg of hydrogen from the electrolyzer output pressure to the pressure level of the storage or tube trailers. The hourly hydrogen demand is bounded by the capacity of tube trailers, as 
\begin{align}
    d_{t} &\leq{D}^\mathrm{max} & \forall~ & t \in \mathcal{T}.
\end{align}

In case the hydrogen production $h_t^{\mathrm{p}}$ of the electrolyzer in hour $t$ is more than demand $d_{t}$ in that hour, the excess is being stored, while in the case of deficit, we discharge the storage. By this, the state of charge of the storage $h^{\text{s}}_{t}$ is defined as 
\begin{align}
    h^{\text{s}}_{t} &=h^{\text{p}}_{t} - d_{t} \label{eq:tank_start} & &t =1,\\
    h^{\text{s}}_{t} &=h^{\text{p}}_{t} - d_t + h^{\text{s}}_{t-1}  &\forall~ & t \in \mathcal{T}\backslash 1, \label{eq:tank_gen} 
\end{align}
which is upper-bounded by the capacity of the storage, i.e., 
\begin{align}
    h^{\text{s}}_{t} &\leq   {H}^\mathrm{max}& \forall~ & t \in \mathcal{T}. \label{eq:tank_cap} 
\end{align}
\end{subequations}

The following set of constraints \eqref{eq:FCRbid} enforces FCR reserve allocation constraints. To clarify the need for these constraints, Figure \ref{fig:ele_fcr} provides an example, where the electrolyzer consumes $p^{\text{e}}_{t}$ in hour $t$, which is the baseline for reserve activation. Recall from \eqref{eq:FCR} that FCR-N is a market with a two-side product, meaning that the electrolyzer might be activated to consume less power (if frequency is below 50 Hz) or more power (if frequency is above 50 Hz). On the contrary, FCR-D up/Down are markets with one-side products, meaning that if FCR-D Up is activated (i.e., frequency is below 49.9 Hz), the electrolzyer must consume less power, whereas if FCR-D Down is activated (i.e., frequency is above 50.1 Hz),  the electrolzyer must consume more power. To operate fully reliably under the worst case wherein frequency drops to 49.5 Hz (threshold defined by the Nordic TSOs), the electrolyzer should be able to respond by the full activation of both FCR-N and FCR-D Up, i.e.,
\begingroup
\allowdisplaybreaks
\begin{subequations} \label{eq:FCRbid}
\begin{align}
    p^{\text{e}}_{t} -  r_{t}^{\text{FCR-N}} -  r_{t}^{\text{FCR-D}\uparrow}  &\geq P^{\text{min}}  z^{\text{on}}_{t} + P^{\text{sb}}  z^{\text{sb}}_{t} & \forall~ & t \in \mathcal{T}, \label{eq:FCR_upper} \
\end{align}
indicating that if fully activated, the electrolyzer's consumption should still not be lower than $P^{\text{min}}$ (if it is on) or $P^{\text{sb}}$ if it is in the standby state. Similarly, for the over-frequency worse case (i.e., 50.5 Hz defined by the Nordic TSOs), we enforce 
\begin{align}
    p^{\rm{e}}_t +  r_{t}^{\text{FCR-N}} +  r_{t}^{\text{FCR-D}\downarrow}  &\leq P^{\text{max}}  z^{\text{on}}_{t} + P^{\text{sb}}  z^{\text{sb}}_{t} & \forall~ & t \in \mathcal{T}, \label{eq:FCR_lower}
\end{align}
stating that by the full activation of both FCR-N and FCR-D Down, the electrolyzer's consumption should not go beyond its capacity $P^{\text{max}}$ (if on) or $P^{\text{sb}}$ (if standby)\footnote{An extension to this work is to make \eqref{eq:FCR_upper} and \eqref{eq:FCR_lower} probabilistic, e.g., via chance constraints, making less conservative reserve allocation decisions, which is outside the scope of this paper. We refer the interested reader to \cite{line}.}.


In addition, we enforce the minimum bid requirement $Q^{\text{FCR}}$, an identical value for all FCR services set by the Nordic TSOs:
\begin{align}
    r_{t}^{\text{FCR-D}\uparrow}  &\geq z_{t}^{\text{FCR-D}\uparrow} Q^{\text{FCR}}   & \forall~ & t \in \mathcal{T}, \label{eq:FCRDup_down_bid} \\
    r_{t}^{\text{FCR-D}\uparrow}  &\leq z_{t}^{\text{FCR-D}\uparrow} \left(  P^{\text{max}} - P^{\text{min}} \right)   & \forall~ & t \in \mathcal{T}, \label{eq:FCRDup_up_bid} \\
        r_{t}^{\text{FCR-D}\downarrow}  &\geq z_{t}^{\text{FCR-D}\downarrow} Q^{\text{FCR}}   & \forall~ & t \in \mathcal{T}, \label{eq:FCRDdown_down_bid} \\
    r_{t}^{\text{FCR-D}\downarrow}  &\leq z_{t}^{\text{FCR-D}\downarrow}  \left(  P^{\text{max}} - P^{\text{min}} \right)   & \forall~ & t \in \mathcal{T}, \label{eq:FCRDdown_up_bid} \\
        r_{t}^{\text{FCR-N}}  &\geq z_{t}^{\text{FCR-N}} Q^{\text{FCR}}   & \forall~ & t \in \mathcal{T}, \label{eq:FCRN_down_bid} \\
r_{t}^{\text{FCR-N}}  &\leq z_{t}^{\text{FCR-N}} \left( \frac{ P^{\text{max}} - P^{\text{min}} }{2}\right)   & \forall~ & t \in \mathcal{T}, \label{eq:FCRN_up_bid}
\end{align}
\end{subequations} 
where binary variables $z_{t}^{\text{FCR-N}}$, $z_{t}^{\text{FCR-D }\uparrow}$, and $z_{t}^{\text{FCR-D }\downarrow}$ ensure that the reserve quantity bid $r^{\text{(.)}}_t$, if takes a non-zero value, is lower bounded by $Q^{\text{FCR}}$. If a binary variable $z_{t}^{\text{(.)}}$ takes a zero value, combination of the corresponding lower and upper bounds enforces $r_{t}^{\text{(.)}}$ to be zero. 
%


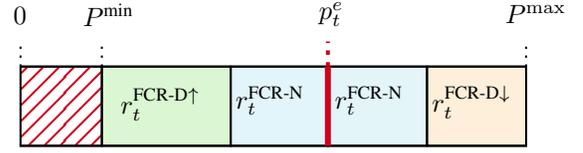
\begin{figure}
    \centering

\usetikzlibrary{patterns}

 
\tikzset{
pattern size/.store in=\mcSize, 
pattern size = 5pt,
pattern thickness/.store in=\mcThickness, 
pattern thickness = 0.3pt,
pattern radius/.store in=\mcRadius, 
pattern radius = 1pt}
\makeatletter
\pgfutil@ifundefined{pgf@pattern@name@_1msdacnmj}{
\pgfdeclarepatternformonly[\mcThickness,\mcSize]{_1msdacnmj}
{\pgfqpoint{0pt}{0pt}}
{\pgfpoint{\mcSize+\mcThickness}{\mcSize+\mcThickness}}
{\pgfpoint{\mcSize}{\mcSize}}
{
\pgfsetcolor{\tikz@pattern@color}
\pgfsetlinewidth{\mcThickness}
\pgfpathmoveto{\pgfqpoint{0pt}{0pt}}
\pgfpathlineto{\pgfpoint{\mcSize+\mcThickness}{\mcSize+\mcThickness}}
\pgfusepath{stroke}
}}
\makeatother
\tikzset{every picture/.style={line width=0.75pt}} 

\begin{tikzpicture}[x=0.75pt,y=0.75pt,yscale=-1,xscale=1]

\draw   (38,92) -- (293,92) -- (293,132) -- (38,132) -- cycle ;
\draw  [pattern=_1msdacnmj,pattern size=6pt,pattern thickness=0.75pt,pattern radius=0pt, pattern color={rgb, 255:red, 208; green, 2; blue, 27}] (38,92) -- (79,92) -- (79,132) -- (38,132) -- cycle ;
\draw  [fill={rgb, 255:red, 74; green, 193; blue, 226 }  ,fill opacity=0.15 ] (193,92) -- (243,92) -- (243,132) -- (193,132) -- cycle ;
\draw  [fill={rgb, 255:red, 74; green, 193; blue, 226 }  ,fill opacity=0.16 ] (144,92) -- (193,92) -- (193,132) -- (144,132) -- cycle ;
\draw  [fill={rgb, 255:red, 71; green, 211; blue, 33 }  ,fill opacity=0.2 ] (79,92) -- (144,92) -- (144,132) -- (79,132) -- cycle ;
\draw  [fill={rgb, 255:red, 245; green, 166; blue, 35 }  ,fill opacity=0.17 ] (243,92) -- (293,92) -- (293,132) -- (243,132) -- cycle ;
\draw [color={rgb, 255:red, 208; green, 2; blue, 27 }  ,draw opacity=1 ][line width=2.25]    (193,92) -- (193,132) ;
\draw  [dash pattern={on 0.84pt off 2.51pt}]  (38,92) -- (38,79.33) ;
\draw  [dash pattern={on 0.84pt off 2.51pt}]  (79,92) -- (79,79.33) ;
\draw  [dash pattern={on 0.84pt off 2.51pt}]  (293,92) -- (293,79.33) ;
\draw [color={rgb, 255:red, 208; green, 2; blue, 27 }  ,draw opacity=1 ][line width=1.5]  [dash pattern={on 1.69pt off 2.76pt}]  (193,92) -- (193,79.33) ;

\draw (68,58.9) node [anchor=north west][inner sep=0.75pt]    {$P^{\text{min}}$};
\draw (281,59.4) node [anchor=north west][inner sep=0.75pt]    {$P^{\max}$};
\draw (187,57.9) node [anchor=north west][inner sep=0.75pt]    {$p_{t}^{e}$};
\draw (87,101.4) node [anchor=north west][inner sep=0.75pt]    {$r_{t}^{\text{FCR-D}\uparrow }$};
\draw (244,100.4) node [anchor=north west][inner sep=0.75pt]    {$r_{t}^{\text{FCR-D} \downarrow }$};
\draw (145,100.4) node [anchor=north west][inner sep=0.75pt]    {$r_{t}^{\text{FCR-N}}$};
\draw (195,100.4) node [anchor=north west][inner sep=0.75pt]    {$r_{t}^{\text{FCR-N}}$};
\draw (33,60.4) node [anchor=north west][inner sep=0.75pt]    {$0$};

\end{tikzpicture}
\caption{An example FCR reserve allocation for an alkaline electrolyzer, consuming power $p_t^{\text{e}}$ in hour $t$ which is the operational baseline.}
\label{fig:ele_fcr}
\figspace
\end{figure}


Within a Hydrogen Purchase Agreement (HPA), the electrolyzer owner might be obliged to supply at least a minimum demand $\rm{HPA}^{\text{min}}$ over a time period, e.g., a day, or a week. For example, let $\mathcal{T}$ indicate the time horizon, then $\mathcal{H}_w \subset \mathcal{T} $ represents the time period $w \in \mathcal{W}$ where the hydrogen demand must be met. The minimum hydrogen demand is enforced by 
%
\begin{align}
    \sum_{t\in \mathcal{H}_w} d_t &\geq \rm{HPA}^{\text{min}} & \forall \ w \in \mathcal{W}.
\end{align}


The non-negativity of variables is enforced by
\begin{subequations} \label{eq:var}
\begin{align}
\begin{split}
    &    d_t, h^{\text{s}}_{t}, h^{\text{p}}_t, p_t^{\text{e}}, \widehat{p}^{\text{e}}_{t, s}, p_{t},  p^{\rm{c}}_{t}, r_{t}^{\text{FCR-D}\uparrow}, r_{t}^{\text{FCR-D}\downarrow}, r_{t}^{\text{FCR-N}} \in \mathbb{R}^{+}  .
   \label{eq:non-neg}
\end{split}
\end{align}
whereas binary variables are
\begin{equation}
    z^{\text{on}}_{t}, z^{\text{sb}}_{t}, z_{t}^{\text{su}}, \widehat{z}_{t, s}, z_{t}^{\text{FCR-D}\uparrow}, z_{t}^{\text{FCR-D}\downarrow}, z_{t}^{\text{FCR-N}} \in \{0, 1\}.
    \label{eq:binary}
\end{equation}
\end{subequations} 
The vector $\bf{x}$ contains all variables listed in \eqref{eq:var}.

%% file: Sections/04_Results.tex
\section{Numerical results and analysis}
\label{results}
We consider a 10-MW alkaline electrolyzer located in DK2, operating at a pressure of 30 bars, increased up to 350 bars by the compressor for storage purposes. For the hydrogen production curve of the electrolyzer, we use five linearization segments. The minimum weekly hydrogen demand is 9,072 kg, which can be met if the electrolyzer operates with 30\% of its capacity all over the week.
All parameters are given in Table~\ref{tab:04_table}. Currently, the minimum bid quantity $Q^{\text{FCR}}$ in the Nordic area is 0.1 MW, which is not a limit for a 10-MW electrolyzer, but it could be for small-scale electrolyzers. 

We solve the proposed optimization problem based on real prices in year 2022. Since we use realized (and not forecasted) prices, this study provides an economic assessment assuming a perfect foresight for year 2022. All source codes and input data are
publicly shared\footnote{Github repository: \url{https://github.com/marco-srtt/electrolyzer_nordic_FCR}}.


\begin{table}[t]
\figspace
    \centering
    \caption{Input data.}
    \resizebox{0.35\textwidth}{!}{%
\begin{tabular}{llllc}
\toprule
\multirow{2}{*}{Tariffs}    & $\lambda^{\text{TSO}}$ & 15.6 & [\euro/MWh] & \cite{energinet-elmarkedet}\\
                              & $\lambda^{\text{DSO}}$ & 5.36 & [\euro/MWh]&\cite{radiuselnet}\\
\midrule
\multirow{5}{*}{Electrolyzer}  
                              & $P^{\text{min}}$ &1.6 & [MW]&\\
                            & $P^{\text{max}}$ &10 & [MW]&\\
                           & $P^{\text{sb}}$  & 0.5 & [MW] & \cite{MATUTE20211449}\\
                           & $K^\mathrm{su}$ & 1,000 & [\euro] &\\
                           & $K^\mathrm{c}$  & 1.67 & [kWh/kg \ce{H2}] &\\
\midrule
\multirow{4}{*}{Hydrogen}  & $\lambda^{{\ce{H2}}}$  & 2 & [\euro/kg] &\\
                            & $\rm{HPA}^\mathrm{min}$  & 9,072  & [kg/week] &\\
                           & ${H}^\mathrm{max}$  & 60,500 & [kg] &\\
                           & ${D}^\mathrm{max}$  & 180  & [kg/hour] &\\
\midrule
\multirow{1}{*}{FCR}       & $Q^{\text{FCR}}$ & 0.1 & [MW] &\cite{energinet2012ancillary}\\
\bottomrule
\end{tabular}}
\label{tab:04_table}
\figspace
\vspace{2mm}
\end{table}

\subsection{Optimal electrolyzer scheduling}
\begin{figure}[t]
\figspace
    \centering
    \includegraphics[clip, trim=0.2cm 0.3cm 0.3cm 0.25cm, width =0.5 \textwidth]{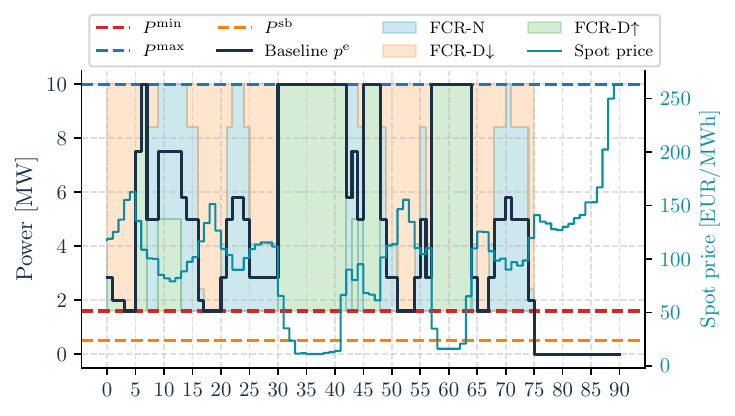}
    \caption{Optimal scheduling of the electrolyzer in FCR-N, FCR-D Up/Down, and spot markets, in an example 90-hour horizon, starting from 22/02/2022.}
\figspace
    \label{fig:ele_op}
\end{figure}

Figure \ref{fig:ele_op} illustrates the electrolyzer scheduling during a sample 90-hour horizon within 2022. We make four observations. 

First, in hours with comparatively low spot prices, e.g., hours 30-40, the electrolyzer operates in its full capacity of 10 MW to maximize hydrogen production. Among FCR services, the electrolyzer sells FCR-D Up reserve only in these hours. The FCR-D Up bid quantity is maximum, which is 8.4 MW, i.e., the capacity of 10 MW minus the minimum operating level of 1.6 MW. 

Second, in hours with comparatively high spot prices, e.g., hours 16-19, the electrolyzer operates in its minimum operating level of 1.6 MW. Among FCR services, the electrolyzer sells FCR-D Down reserve only in these hours. Again,  the electrolyzer submits its maximum FCR-D Down bid quantity, which is 8.4 MW (i.e., 10 MW minus 1.6 MW).

Third, in hours with extremely high spot prices, e.g., hours 75-90, the electrolyzer switches off, producing neither hydrogen nor any FCR services. Indeed, this might be affected if the minimum weekly hydrogen demand is higher, eventually reducing the profit. 

Fourth, in hours with intermediate spot prices, the electrolyzer operates at partial loading between its minimum level of 1.6 MW and the capacity of 10 MW, and produces also FCR-N services. For example, in hours 22-24, among FCR services, it only produces FCR-N. For that, the electrolyzer consumes 5.8 MW to be able to sell 4.2 MW reserve in the FCR-N auction, such that in extreme cases of full activation, the consumption level either drops to the minimum or increases to the maximum level. There are also hours that the electrolyzer sells multiple FCR services, e.g., FCR-N and FCR-D Up in hours 9-13, or FCR-N and FCR-D Down in hours 25-30.  

Note that the minimum weekly demand of 9,072 kg hydrogen is met in the reserve contraction stage. In the activation stage it might be violated, although it is unlikely as already explained in footnote 1. However, one may develop a real-time policy to track meeting the weekly demand, which is outside the scope of this paper.

\begin{figure}[b]
\figspace
    \centering
    \includegraphics[clip, trim=1.1cm 2.8cm 1.1cm 2.8cm, width =0.5 \textwidth]{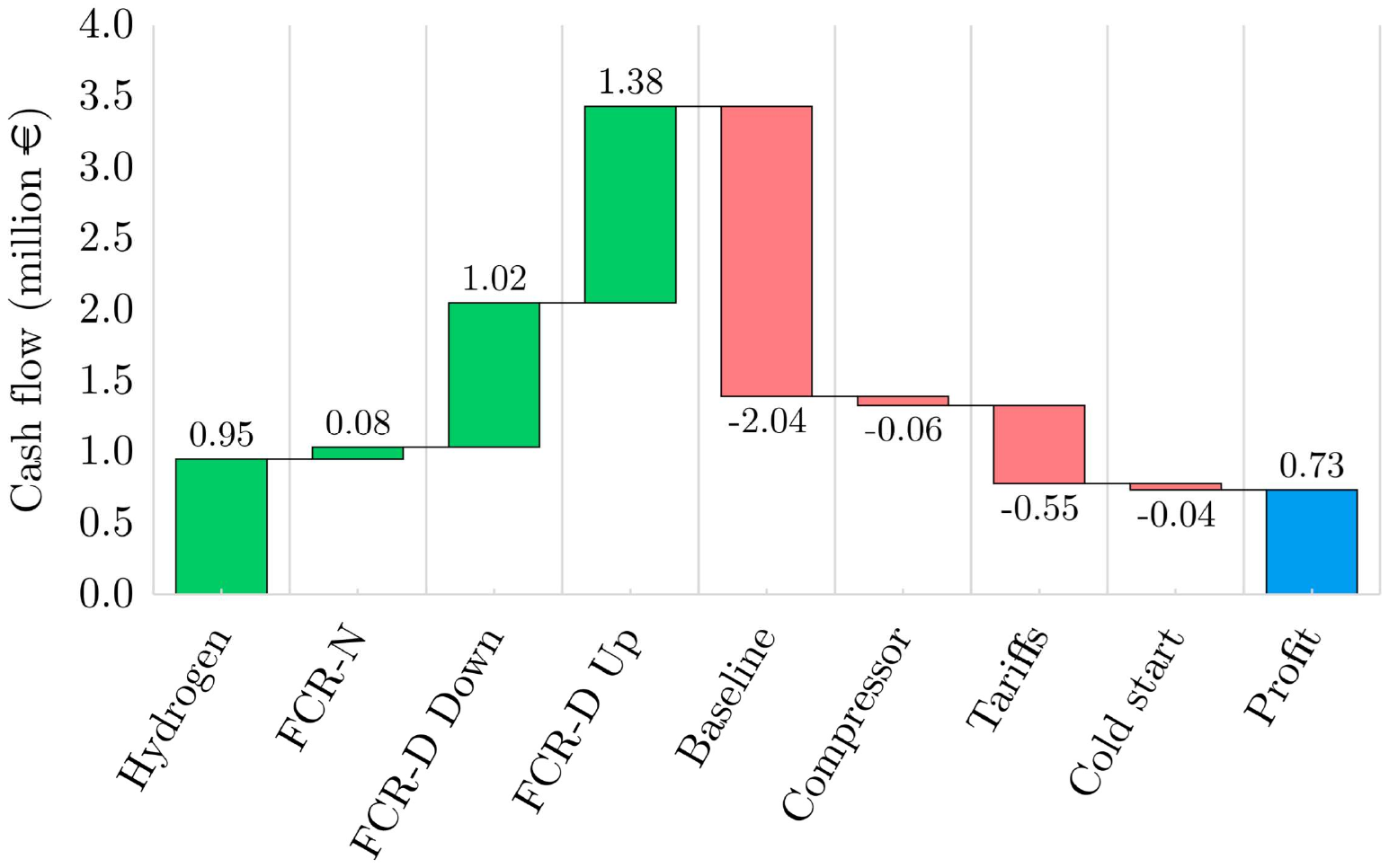}
    \caption{Cash flow for a 10-MW alkaline electrolyzer, participating in the Nordic FCR markets in 2022. Minimum weekly hydrogen demand is 9,072 kg, equivalent to 30\% of electrolyzer's capacity. The hydrogen price is \euro2/kg.}
    \label{fig:waterfall}
    \figspace
\end{figure}
\subsection{Economic assessment}
Figure \ref{fig:waterfall} shows the yearly profit of the electrolyzer in 2022, which is 0.73 million \euro, as well as the distribution of yearly revenues and expenses. The activation payments are excluded, but as it was mentioned earlier, FCR services are not energy-intensive overall, and thereby the activation payments are expected to be negligible \cite{thingvad2022economic}. 

The total annual revenue is 3.43 million \euro, for which the contributions of selling hydrogen, FCR-N, FCR-D Down, and FCR-D Up are  28\%, 2\%, 30\%, and 40\%, respectively. This implies that the electrolyzer earns 72\% of its total revenue from FCR auctions, which is significant. Indeed, these results could be sensitive to the hydrogen price of \euro2/kg and the minimum weekly hydrogen demand of 30\%. Therefore, we will conduct a sensitivity analysis in the next section. 
The total expenses over the year 2022 is 2.69 million \euro, 76\% of which corresponds to the power consumption of the electrolyzer (baseline). Tariffs cause 20\% of total expenses, while the remaining 4\%  is incurred by the consumption of the compressor and the start-up cost of electrolyzer (44 start-ups over 2022, each costing 1,000 \euro). 

\subsection{Sensitivity analysis}
Recall we have assumed the minimum weekly hydrogen demand $\rm{HPA}^\mathrm{min}$ is met if the electrolyzer operates at 30\% of its capacity all over the week, whereas the hydrogen price $\lambda^{{\ce{H2}}}$ is \euro2/kg. We conduct a sensitivity analysis for the annual profit of the electrolyzer with respect to these two parameters. We vary $\rm{HPA}^\mathrm{min}$ from 0\% to 50\%, and $\lambda^{{\ce{H2}}}$ from \euro1/kg to \euro5/kg. We conduct this analysis for two cases: (\textit{i}) the electrolyzer offers FCR services along with the hydrogen production, and (\textit{ii}) the electrolyzer produces hydrogen only. 



\begin{figure}[]
    \centering
    \includegraphics[width =0.5 \textwidth]{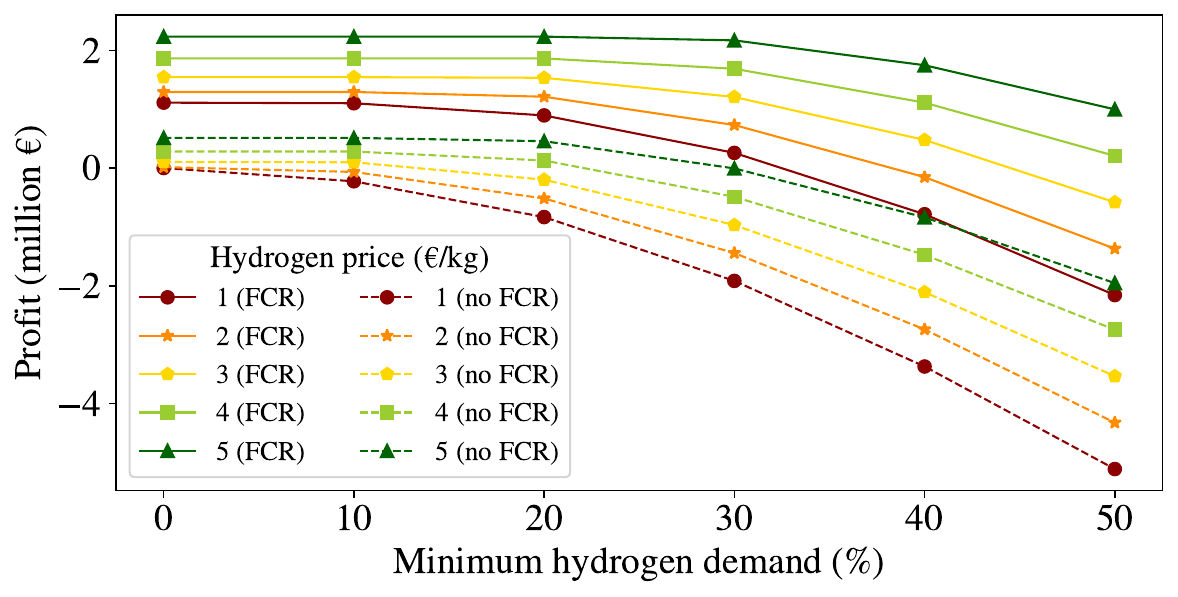}
    \caption{Sensitivity of the annual profit of the electrolyzer (with and without selling FCR services) with respect to the minimum weekly hydrogen demand (\%) and the hydrogen price (\euro/kg).}
    \label{fig:sens}
    \figspace
\end{figure}

The results are depicted in Figure \ref{fig:sens}. As expected, profit declines by increasing $\rm{HPA}^\mathrm{min}$, as the electrolyzer is obliged to operate during non-profitable hours. In extreme cases, the annual profit is negative. The economic value of FCR services becomes even more remarkable when $\rm{HPA}^\mathrm{min}$ increases. Finally, higher hydrogen prices increase the profit.



%% file: Sections/05_Conclusion.tex
\section{Conclusion}
\label{conclusion}

This paper develops a mixed-integer linear model for optimal scheduling of an electrolyzer, purchasing power from the spot market and selling hydrogen as well as FCR-N and FCR-D Up/Down services in the Nordic synchronous region. For a case study based on realized spot and FCR prices in 2022, we found out FCR services can significantly increase the annual profit of the electrolyzer. For a case with the fixed hydrogen price of \euro2/kg and the minimum weekly hydrogen demand of 30\%, the electrolyzer earns the annual profit of 0.73 million \euro \ with a significant contribution from FCR markets (72\%), particularly from FCR-D Up/Down markets. However, this is an analysis with perfect foresight into prices, thereby the true contribution of FCR markets with imperfect foresight might be different. 
The capital cost of an alkaline electrolyzer alone (without auxiliary assets) varies with its scale and depends on the manufacturer, but overall it is approximately around 1 million \euro \ per MW \cite{grigoriev2020current}. It  looks the annual profit of 0.73 million \euro, earned mostly from FCR auctions, is still insufficient to recover the investment cost, but it requires an in-depth analysis, which is left for the future work. This may call for additional regulatory supportive actions to make green hydrogen cost competitive.

%% file: conference_101719.bbl
\begin{thebibliography}{10}

\bibitem{PtX_strategy_DK}
{Danish Ministry of Climate, Energy and Utilities}, ``The goverment's strategy for power-to-x,'' 2021.
\newblock \url{https://ens.dk/sites/ens.dk/files/ptx/strategy_ptx.pdf}.

\bibitem{GUERRA20192425}
O.~J. Guerra {\em et~al.}, ``Cost competitiveness of electrolytic hydrogen,'' {\em Joule}, vol.~3, no.~10, pp.~2425--2443, 2019.

\bibitem{newbusinessmodels}
P.~Larscheid {\em et~al.}, ``Potential of new business models for grid integrated water electrolysis,'' {\em Renew. Energy}, vol.~125, pp.~599--608, 2018.

\bibitem{VARELA20219303}
C.~Varela, M.~Mostafa, and E.~Zondervan, ``Modeling alkaline water electrolysis for power-to-x applications: A scheduling approach,'' {\em Int. J. Hydrog. Energy}, vol.~46, no.~14, pp.~9303--9313, 2021.

\bibitem{BUTTLER20182440}
A.~Buttler and H.~Spliethoff, ``Current status of water electrolysis for energy storage, grid balancing and sector coupling via power-to-gas and power-to-liquids: A review,'' {\em Renew. Sustain. Energy Rev.}, vol.~82, pp.~2440--2454, 2018.

\bibitem{Mancarella}
M.~G. Dozein, A.~Jalali, and P.~Mancarella, ``Fast frequency response from utility-scale hydrogen electrolyzers,'' {\em IEEE Trans. Sustain. Energy}, vol.~12, no.~3, pp.~1707--1717, 2021.

\bibitem{Tuinema}
B.~W. Tuinema {\em et~al.}, ``Modelling of large-sized electrolysers for real-time simulation and study of the possibility of frequency support by electrolysers,'' {\em IET Gener. Transm. Distrib}, vol.~14, no.~10, pp.~1985--1992, 2020.

\bibitem{ALLIDIERES20199690}
L.~Allidières {\em et~al.}, ``On the ability of {PEM} water electrolysers to provide power grid services,'' {\em Int. J. Hydrog. Energy}, vol.~44, no.~20, pp.~0360--3199, 2019.

\bibitem{France}
B.~Guinot {\em et~al.}, ``Profitability of an electrolysis based hydrogen production plant providing grid balancing services,'' {\em Int. J. Hydrog. Energy}, vol.~40, no.~29, pp.~8778--8787, 2015.

\bibitem{Germany_offshore}
M.~Scolaro and N.~Kittner, ``Optimizing hybrid offshore wind farms for cost-competitive hydrogen production in {G}ermany,'' {\em Int. J. Hydrog. Energy}, vol.~47, no.~10, pp.~6478--6493, 2022.

\bibitem{ZHENG2023}
Y.~Zheng {\em et~al.}, ``Economic evaluation of a power-to-hydrogen system providing frequency regulation reserves: {A} case study of {D}enmark,'' {\em Int. J. Hydrog. Energy}, 2023.

\bibitem{thingvad2022economic}
A.~Thingvad {\em et~al.}, ``Economic value of multi-market bidding in {N}ordic frequency markets,'' in {\em 2022 International Conference on Renewable Energies and Smart Technologies}, pp.~1--5, 2022.

\bibitem{energinet2012ancillary}
Energinet, ``Ancillary services to be delivered in {D}enmark: Tender conditions,'' 2023.
\newblock \url{https://en.energinet.dk/Electricity/Ancillary-Services/Tender-conditions-for-ancillary-services/}.

\bibitem{fingrid}
{Fingrid}, ``{Frequency - historical data}.'' \url{https://data.fingrid.fi/en/dataset/frequency-historical-data}.

\bibitem{baumhof2023optimization}
M.~T. Baumhof {\em et~al.}, ``Optimization of hybrid power plants: When is a detailed electrolyzer model necessary?,'' in {\em 2023 IEEE Belgrade PowerTech}, pp.~1--10, 2023.

\bibitem{line}
A.~Porras {\em et~al.}, ``Integrating automatic and manual reserves in optimal power flow via chance constraints,'' 2023.
\newblock \url{https://arxiv.org/abs/2303.05412}.

\bibitem{energinet-elmarkedet}
{Energinet}, ``{Aktuelle Tariffer}.'' \url{https://energinet.dk/el/elmarkedet/tariffer/aktuelle-tariffer/}.

\bibitem{radiuselnet}
{Radius Elnet}, ``{Tariffer og Netabonnement}.'' \url{https://radiuselnet.dk/elnetkunder/tariffer-og-netabonnement/}.

\bibitem{MATUTE20211449}
G.~Matute {\em et~al.}, ``Multi-state techno-economic model for optimal dispatch of grid connected hydrogen electrolysis systems operating under dynamic conditions,'' {\em Int. J. Hydrog. Energy}, vol.~46, no.~2, pp.~1449--1460, 2021.

\bibitem{grigoriev2020current}
S.~Grigoriev {\em et~al.}, ``Current status, research trends, and challenges in water electrolysis science and technology,'' {\em Int. J. Hydrog. Energy}, vol.~45, no.~49, pp.~26036--26058, 2020.

\end{thebibliography}
